\input amstex

\loadeufm
\loadmsbm
\loadeufm

\documentstyle{amsppt}
\input amstex
\catcode `\@=11
\def\logo@{}
\catcode `\@=12
\magnification \magstep1
\NoRunningHeads
\NoBlackBoxes
\TagsOnLeft

\def \={\ = \ }
\def \+{\ +\ }
\def \-{\ - \ }

\def \b|{\big |}

\def \g1{\Gamma_1}

\def \nfp{\demo\nofrills{Proof:\usualspace\usualspace }}

\def\rarr#1#2{\smash{\mathop{\hbox to .5in{\rightarrowfill}}
 	 \limits^{\scriptstyle#1}_{\scriptstyle#2}}}

\def\larr#1#2{\smash{\mathop{\hbox to .5in{\leftarrowfill}}
	  \limits^{\scriptstyle#1}_{\scriptstyle#2}}}

\def\swarr#1#2 {\llap{$\scriptstyle #1$}  \swarrow
  	\vcenter to .5in{}\rlap{$\scriptstyle #2$}}

\topmatter
\title Free and Hyperbolic Groups are Not Equational 
\endtitle
\author
\centerline{ 
Z. Sela${}^{1,2}$}
\endauthor
\footnote""{${}^1$Hebrew University, Jerusalem 91904, Israel.}
\footnote""{${}^2$Partially supported by an Israel academy of sciences fellowship.}
\abstract\nofrills{}
In [Se5] we proved that free and torsion-free hyperbolic groups are stable. In this note we give an example of a definable
set in each of these groups that is not in the Boolean algebra of equational sets. Hence, the theories of free and torsion-free
hyperbolic groups are not equational in the sense of G. Srour.
\endabstract
\endtopmatter

\document

\baselineskip 12pt

In [Se5] we proved that free groups and torsion-free hyperbolic 
groups are stable. We started by proving that a subset
of the collection of definable sets in these groups, that are called minimal rank definable sets, are
in the Boolean algebra of equational sets, and that general varieties and Diophantine sets in these groups
are equational.
Then we used the equationality of Diophantine sets and Duo limit groups and their properties (that are
presented in section 3 of [Se5]) to prove that free
and torsion-free hyperbolic groups are stable.

In this note, we show that our results for minimal rank definable sets are false in general. We give an
example of a definable set in each free or non-elementary torsion-free hyperbolic group that is not in the Boolean 
algebra of equational sets.  

Recall that equational sets and theories were defined by Gabriel Srour (see [Pi-Sr]).
A definable set  $D(p,q)$ is called $equational$ if there exists
a constant $N_{D}$,
so that  for every sequence of values $\{q_i\}_{i=1}^{m}$, for which the sequence of sets
that corresponds to the intersections:
$ \{ \, \cap_{i=1}^{j} D(p,q_i) \, \}_{j=1}^m$ is a strictly decreasing sequence, satisfies:
$m \leq N_{D}$.  

We note that the question of the existence
of a theory 
which is stable but not equational was raised by Pillay and Srour [Pi-Sr], and such 
examples were constructed by Hrushovski and Srour [Hr-Sr] and Baudisch and Pillay [Ba-Pi].

\vglue 1pc
\proclaim{Theorem 1} The elementary theory of a non-abelian free group is
not equational.
\endproclaim

\nfp Let $F_k$ be a non-abelian free group, and let
$NE(p,q)$ be the existential set:
$$ NE(p,q) \ = \ \exists x_1,x_2 \ \ qp=x_1^{10}x_2^{-9} \, \wedge \, [x_1,x_2] \neq 1.$$
In section 3 of [Se5] we introduced duo limit groups (see definition 3.1 in [Se5]).
With the set $NE(p,q)$ we associate the following natural duo limit group:
$$Duo \, = \, <d_1> *_{<d_0>} \, <d_2> \ = \ <t,u>*_{<u>} \, <u,s> \, ; \ q=s^{10}u \, ; \ p = u^{-1}t^{-9}.$$
Note that if we fix a value of the variable $u$, i.e., we fix a rectangle that is associated
with the duo limit group $Duo$, then for every (fixed) non-trivial value of the variable $s$, a
couple $(q(s),p(t))$ is in the set $NE(p,q)$ for every generic value of the variable $t$,
and correspondingly for every (fixed)  non-trivial value of $t$ and a generic value of $s$. 

\vglue 1pc
\proclaim{Lemma 2} The set $NE(p,q)$ is not equational.
\endproclaim

\nfp If we fix a value $u_0$ for the variable $u$, and a non-trivial value $s_0$ for the
variable $s$, so that $q(s_0)=s_0^{10}u_0$, then for all values of the variable $t$ for which
$[t,s_0] \neq 1$, the couple $(q(s_0),p(t)) \in NE(p,q)$ where $p(t)=u_0^{-1}t^{-9}$. A work of Lyndon and
Schutzenberger [Ly-Sch], shows that for $m,n,r \geq 2$, the solutions
of the equation $x^my^nz^r=1$ in a free group, generate a cyclic subgroup. Hence, for all
but at most 2 values $t_0$ of the variable $t$, for which $[s_0,t_0]=1$, $(q(s_0),q(t_0)) \notin
NE(p,q)$. Hence, if we consider a sequence, $s_1,s_2,\ldots$, where $[s_{i_1},s_{i_2}] \neq 1$
for $i_1<i_2$, then the sequence of intersections:
$\{ \cap_{i=1}^j \, NE(p,q(s_i)) \}_{j=1}^{\infty}$ is a strictly decreasing sequence. Therefore,
$NE(p,q)$ is not equational.

\line{\hss$\qed$}

Since a finite union and a finite intersection of equational sets are equational (see corollary 2.27 in [OH]),
lemma 2 implies that the set $NE(p,q)$ is not
a finite  union nor a finite intersection of equational sets. 

\vglue 1pc
\proclaim{Lemma 3} The set $NE(p,q)$ is not co-equational (i.e., it is not the complement of an
equational set).
\endproclaim

\nfp The (infinite) union of the sets $NE(p,q_i)$, for all possible values of $q_i$, is the entire 
coefficient group $F_k$. Every finite subunion of such sets is a proper subset of $F_k$ (since the set $NE(p,q)$ is not
of generic type or is negligible in the sense of M. Bestvina and M. Feighn [Be-Fe]). Hence, 
$NE(p,q)$ can not be co-equational.

\line{\hss$\qed$}
 
To prove that the theory of a free group is not equational, we need to show that the set
$NE(p,q)$ is  not a union nor an intersection of an equational and a co-equational sets.

\vglue 1pc
\proclaim{Lemma 4} The set $NE(p,q)$ is not a union of an equational and a co-equational sets.
\endproclaim

\nfp Suppose that $NE(p,q) \, = \, E_1 \, \cup \, (E_2)^c$ where $E_1,E_2$ are equational sets.
Let $N_1,N_2$ be the bounds on the lengths of 
sequences with strictly decreasing intersections: $\{E_1(p,q_i)\}$ and $\{E_2(p,q_j)\}$.
Let $N= \max (N_1,N_2)$.
We look at a sequence of values: $\{q^i_j={(s^i_j)}^{10}u_i\}_{i,j=1}^{N+1}$, 
and
$\{p^i_j=u_i^{-1}{(s^i_j)}^{-9}\}_{i,j=1}^{N+1}$, 
where $\{s^i_j\}$
and $\{u_i\}$ are distinct elements in a given test sequence 
(i.e., they satisfy the combinatorial properties that are listed in definition 1.20 in [Se2] for a test sequence. 
Elements in a test sequence should be 
considered "generic" elements, that as a set satisfy a very small cancellation condition).

Since the elements $\{u_i\}$ and the elements $\{s^i_j\}$ are chosen to be generic (part of a test sequence), a pair,
$(p^i_j,q^{i'}_{j'}) \in NE(p,q)$ if and only if $i=i'$ and $j \neq j'$. Hence, for a given index $i$, $1 \leq i \leq N+1$,
and every index $j$, $1 \leq j \leq N+1$, $E_2(p,q^i_j)$ contains all the elements $p^{i'}_j$, $i \neq i'$, 
$1 \leq i,i',j < N+1$.

Our choice of $N$ guarantees that a  sequence of sets of the form, $\{E_2(p,q^i_j)\}$, 
with strictly decreasing intersections,
 has length
bounded by $N$. Therefore, 
there must exist some index $i_0$, $1 \leq  i_0 \leq N+1$, for which for every
$1 \leq j,j' \leq N+1$:
$p^{i_0}_j \in E_2(p,q^{i_0}_{j'})$.

%\noindent
Hence, for this index $i_0$, $p^{i_0}_j \in E_1(p,q^{i_0}_{j'})$, if and only if $p^{i_0}_j \in NE(p,q^{i_0}_{j'})$,
$1 \leq j,j' \leq N+1$.
$p^{i_0}_j \in NE(p,q^{i_0}_{j'})$ if and only if $j \neq j'$. Therefore,
the sequence of intersections:
$\{ \cap_{j=1}^{m} \, E_1(p,q^{i_0}_j) \}_{m=1}^{N+1}$ is a strictly decreasing sequence, a contradiction
to the assumption that $N$ is the global bound on the length of such a strictly decreasing sequence.

\line{\hss$\qed$}

Using a somewhat similar argument one can show that $NE(p,q)$ is not the intersection of an 
equational and a co-equational sets.

\vglue 1pc
\proclaim{Lemma 5} The set $NE(p,q)$ is not an intersection between an equational and a co-equational sets.
\endproclaim

\nfp Suppose that $NE(p,q) \, = \, E_1 \, \cap \, (E_2)^c$ where $E_1,E_2$ are equational sets.
Let $N_1,N_2$ be the bounds on the lengths of 
strictly decreasing sequences: $\{E_1(p,q_i)\}$ and $\{E_2(p,q_j)\}$,
and let $N= \max (N_1,N_2)$.
We look at a sequence of values: $\{s_i=(v_i)^{10}(w_i)^{-9}\}_{i=1}^{N+1}$,
$\{q^i_j=s_i^{(10^j)}\}_{i,j=1}^{N+1}$, $\{p^i_j=s_i^{(1-10^j)}\}$, 
where $\{v_i\}$
and $\{w_i\}$ are distinct elements in a given test sequence ("generic" elements). 

Since the elements $\{v_i\}$ and  $\{w_i\}$ are chosen to be generic (part of a test sequence), a pair,
$(p^i_j,q^{i'}_{j'}) \in NE(p,q)$, if and only if $i \neq i'$ or $i=i'$ and $j = j'$. Hence, for a given index $i$, 
$1 \leq i \leq N+1$,
and every index $j$, $1 \leq j \leq N+1$, $E_1(p,q^i_j)$ contains the element $p^i_j$, and all the elements 
$p^{i'}_{j'}$, $i \neq i'$, 
$1 \leq i',j' < N+1$.

Our choice of $N$ guarantees that a  sequence of sets of the form, $\{E_1(p,q^i_j)\}$, 
with strictly decreasing intersections,
 has length
bounded by $N$. Therefore, 
there must exist some index $i_0$, $1 \leq  i_0 \leq N+1$, for which for every
$1 \leq j,j' \leq N+1$:
$p^{i_0}_j \in E_1(p,q^{i_0}_{j'})$.

For this index $i_0$, $p^{i_0}_j \in (E_2)^c(p,q^{i_0}_{j'})$, if and only if $p^{i_0}_j \in NE(p,q^{i_0}_{j'})$,
$1 \leq j,j' \leq N+1$.
$p^{i_0}_j \in NE(p,q^{i_0}_{j'})$ if and only if $j = j'$. Hence,
$p^{i_0}_j \in E_2(p,q^{i_0}_{j'})$ if and only if $j \neq j'$. Therefore,
the sequence of intersections:
$\{ \cap_{j=1}^{m} \, E_1(p,q^{i_0}_j) \}_{m=1}^{N+1}$ is a strictly decreasing sequence, a contradiction
to the assumption that $N$ is the global bound on the length of such a strictly decreasing sequence.

\line{\hss$\qed$}
 
So far we proved that the set $NE(p,q)$ is not an intersection nor a union of an equational and
a co-equational sets. Since a finite union and a finite intersection of equational sets are equational,
to prove that the set $NE(p,q)$ is not in the Boolean algebra of equational sets it is enough to
prove that $NE(p,q)$ is not a finite union of intersections between an  equational and a co-equational sets. 
Hence, to prove that the set $NE(p,q)$ is not in the Boolean algebra of equational sets, for every positive 
integer $\ell$, we need to find a collection of elements in the definable set $NE(p,q)$, so that this collection
of elements can not be contained in a union of $\ell$ intersections of an equational set and a co-equational
set.

Unlike the previous lemmas, to prove that we (indirectly) 
use the structure of a definable set over a free group that is the outcome of the quantifier elimination
procedure [Se3]. This (geometric) structure of definable sets
was used in proving the stability of the theory (see theorems 1.9 and 5.1 in [Se5]).
Technically,  to prove theorem 1 (proposition 6 and theorem 7) we use $envelopes$ of definable sets over a free group, 
as they appeared and proved to exist in section 1 in [Se6]. For presentation purposes we start by proving that the set $NE(p,q)$
is not the union of two intersections of equational and co-equational sets (proposition 6), and then generalize the argument to
prove that $NE(p,q)$ is not a finite union of such intersections (theorem 7).  

\vglue 1pc
\proclaim{Proposition 6} The set $NE(p,q)$ is not  a  union of two intersections between an equational and 
a co-equational sets.
\endproclaim

\nfp Our approach to prove the proposition combines the argument that was used to prove lemma 5 with the
envelope of a definable set and its properties as they are presented and   proved to exist 
in theorem 1.3 in [Se6].

Suppose that: $$NE(p,q) \, = \, (E^1_1 \, \cap \, (E^1_2)^c) \, \cup \,  
(E^{2}_1 \, \cap \, (E^{2}_2)^c).$$ 
Let $N$ be the maximum in the finite set of equationality constants of the sets, $E^1_i$, $E^2_i$, 
and $E^1_i \cap E^2_i$, $i=1,2$.

We modify the argument that was used to prove lemma 5, and the definition of the elements $p^i_j$ and $q^i_j$ 
in the proof of that lemma.
We set: 
$$q \, = \, (((t_1)^{10}(s_1)^{-9})^{10m_1} \, ((t_2)^{10}(s_2)^{-9})^{-9m_2})^{10f_1} \, u$$ 
$$p \, = \, u^{-1} \, (((t_3)^{10}(s_3)^{-9})^{10m_3} \, ((t_4)^{10}(s_4)^{-9})^{-9m_4})^{-9f_2}.$$ 

We look at the collection of all the test sequences of values for the variables, 
$\{t_j,s_j\}_{j=1}^{4}$ and $u$,
(recall that elements in a test sequence need to satisfy the properties that are listed in definition 1.20 in [Se2],
they satisfy a very small cancellation assumption, and should be 
considered "generic" elements). For all the values in such test sequences, and arbitrary positive integer values
of the elements, $\{m_j\}_{j=1}^4$ and $f_1,f_2$, the corresponding values of the variables $(p,q)$ 
are in the set $NE(p,q)$.

By our assumptions, for each $(p_0,q_0) \in NE(p,q)$, 
$(p_0,q_0) \in \, E^i_1 \cap (E^i_2)^c$ for $i=1$ or $i=2$, hence, in particular, $(p_0,q_0) \in E^i_1$
for $i=1$ or $i=2$. For $i=1,2$,
 we look at the collection of all the  test sequences of elements, 
$\{t_j,s_j\}_{j=1}^4$ and $u$,
and  sequences of positive integers as values for the elements,
$\{m_j\}_{j=1}^4$ and $f_1,f_2$,
for which each of the sequences of values of each of the elements, 
$\{m_j\}_{j=1}^4$ and $f_1,f_2$,
 converges to infinity,
and the corresponding sequence of values of the pair $(p,q)$ are all in the set $E^i_1$.

By definition,
the group that is generated by the elements, 
$\{t_j,s_j\}_{j=1}^4$ and $u,p,q$,  together with elements that are (formally) added for the powers:
$((t_j)^{10}(s_j)^{-9})^{10m_j}$, $j=1,\ldots,4$,  
(each such formal power commutes with the corresponding element: 
$(t_j)^{10}(s_j)^{-9}$, $j=1,\ldots,4$),  
is a completion (i.e., an $\omega$-residually free tower - see section 6 in [Se1] and definition 1.12 in [Se2]), 
with a bottom level which is a free group of rank 9 (that is generated
by the elements,
$t_j,s_j$, $j=1,\ldots,4$ and $u$,  
and two upper levels, where at the first level above the bottom level, there are 4 abelian vertex groups of rank 2 that are
connected to a distinguished vertex,
and in the top level there are additional two abelian vertex groups of rank 2 that are connected to a distinguished vertex.
We denote this
completion, $Comp$.

The subgroup  $Comp$ is a completion, and  generic points (test sequences) in its
associated variety project to values of the pair, $(p,q)$, that are in $NE(p,q)$, and hence are in $E^i_1(p,q)$  
for $i=1$ or $i=2$. Therefore, by the construction of envelopes of definable sets, as they
are presented in theorem 1.3 in  [Se6], with the sets, $E^i_1(p,q)$, $i=1,2$, 
there are associated envelopes,
that are closures of the completion, $Comp$.

Hence, with  the  sets $E^1_1,E^2_1$, there are associated closures of
$Comp$, that we denote $Cl_1,\ldots,Cl_d$, that form the envelopes of the sets $E^1_1,E^2_1$, with
respect to generic points (test sequences) of the completion, $Comp$. 
Furthermore, with each closure, $Cl_g$, there exists a finite (possibly empty)
closures of it, $Cl_{g,1},\ldots,Cl_{g,f_g}$. By the properties of envelopes (see theorem 1.3 in [Se6]),
the set of closures $Cl_1,\ldots,Cl_d$ and their associated closures,$Cl_{1,1},\ldots,Cl_{g,f_g}$, specifies precisely 
what test sequences of the completion, $Comp$, restrict to values of the variables $(p,q)$ that are in the sets,
$E^1_1$ and $E^2_1$ (and what test sequences restrict to values of $(p,q)$ that are
in the complements of each of these sets). 

\medskip
We now look at a specific  test sequence of values for the variables, 
$t_j,s_j$, $j=1,\ldots,4$ and $u$,  
together with  a specific sequence of positive integers as values for the elements,
$\{m_j\}_{j=1}^4$ and $f_1,f_2$,
for which: $m_1=m_3$, $m_2=m_4$, $10m_1-9m_2=1$, and $10f_2-9f_1=1$. 
 In the sequel we will denote this  specific test
sequence, $\{t_j(n),s_j(n)\}$, $j=1,\ldots,4$, $u(n)$, and $m_j(n)$, $j=1,\ldots,4$ and $f_1(n),f_2(n)$.

With the collection of all the test sequences of the completion, $Comp$, we have associated finitely
many envelope closures, $Cl_1,\ldots,Cl_d$, and with each closure $Cl_g$ we have further associated
finally many (possibly no) closures of it,  $Cl_{g,1},\ldots,Cl_{g,f(g)}$. With each of these closures 
there is an associated finite index subgroup of some fixed f.g.\ free abelian group, so we may look
at the intersection of all these finite index subgroups. Hence, with each closure there are
finitely many associated cosets of this finite index subgroup. By passing to a subsequence we may assume that
all the elements in the specific test sequence that we chose belong to a fixed coset of that finite index subgroup.

Now, for each index $n$ (of the specific test sequence that we chose), we look at the following sequence
of values for $1 \leq g,h \leq N+1$ (that is similar to the one that we used in proving lemma 5):
$$q^g_h(n) \, = \, (((t_1(n+g))^{10}(s_1(n+g))^{-9})^{10m_1(n)} \, ((t_2(n+g))^{10}(s_2(n+g))^{-9})^{-9m_2(n)})^{10f_1(n+h)} \, u(n)$$ 
$$p^g_h(n) \, = \, u(n)^{-1} \, (((t_1(n+g))^{10}(s_1(n+g))^{-9})^{10m_1(n)} \, ((t_2(n+g))^{10}(s_2(n+g))^{-9})^{-9m_2(n)})^{-9f_2(n+h)}.$$ 

By the way the values of the  sequence is defined,
and like in lemma 5, for every index $n$:
$(p^g_h(n),q^{g'}_{h'}(n)) \in NE(p,q)$, if and only if $g \neq g'$ or $g=g'$ and $h = h'$ 
($1 \leq g,g',h,h' \leq N+1$).

Since we have chosen the positive integers, $m_j(n)$, $j=1,\ldots,4$, and $f_1(n),f_2(n)$, to be from a fixed 
class of the finite
index subgroup that is associated with the closures, $Cl_1,\ldots,cl_d$, and
$Cl_{1,1},\ldots,Cl_{d,f(d)}$, for every index $n$, and every $g \neq g'$, $1 \leq h,h',g,g' \leq N+1$, 
precisely one
of the following two possibilities is valid:
\roster 
\item"{(1)}" $(p^g_h(n),q^{g'}_{h'}(n)) \in E^1_1 \cap E^2_1$. 

\item"{(2)}" there exists a fixed index $i$, (wlog  we may assume $i=1$), 
so that: 
$(p^g_h(n),q^{g'}_{h'}(n)) \in E^1_1$ and
$(p^g_h(n),q^{g'}_{h'}(n)) \notin E^2_1$.
\endroster

If possibility (1) is valid, then we apply the same argument that was used to prove lemma 5, 
to the 
equational set $E^1_1 \cap E^2_1$, and the co-equational set, $(E^1_2 \cap E^2_2)^c$, 
and obtain a contradiction. Hence, we may assume that possibility (2) is valid. 

If possibility (2) holds, then by the equationality of the set $E^1_1$ and the argument that was used to
prove lemma 5,  there  exists an index $g_0$, $1 \leq g_0 \leq N+1$, for which for every $h,h'$, $1 \leq h,h' \leq N+1$, 
$(p^{g_0}_h(n),q^{g_0}_{h'}(n)) \in E^1_1$. Since
$(p^{g_0}_h(n),q^{g_0}_{h'}(n)) \in NE(p,q)$ if and only if $h = h'$, the equationality of the set
$E^1_2$ and the argument that was used in proving lemma 5, imply that there exists an index $h_0$, 
$1 \leq h_0 \leq N+1$,  for which for every index $h'$, $1 \leq h' \leq N+1$,
$(p^{g_0}_{h'}(n),q^{g_0}_{h_0}(n)) \in E^1_2$. In particular,  
$(p^{g_0}_{h_0}(n),q^{g_0}_{h_0}(n)) \in E^1_2$. Since   
$(p^{g_0}_{h_0}(n),q^{g_0}_{h_0}(n)) \in NE(p,q) \cap E^1_1 \cap E^1_2$  and:
$$NE(p,q) \, = \, (E^1_1 \, \cap \, (E^1_2)^c) \, \cup \,  
(E^{2}_1 \, \cap \, (E^{2}_2)^c)$$ 
it finally follows that
$(p^{g_0}_{h_0}(n),q^{g_0}_{h_0}(n)) \in E^1_1 \cap E^2_1$.

Now, note that the pairs,  
$(p^{g_0}_{h_0}(n),q^{g_0}_{h_0}(n))$,  can be written as: 
$$q^{g_0}_{h_0}(n) \, 
 = \, (((t_1(n+g_0))^{10}(s_1(n+g_0))^{-9})^{10m_1(n)} \, ((t_2(n+g_0))^{10}(s_2(n+g_0))^{-9})^{-9m_2(n)})^{10f_1(n+h_0)} \, u(n) \, =$$
$$= \, ((t_1(n+g_0))^{10}(s_1(n+g_0))^{-9})^{10m_1(n)} \,  \hat u(n)$$
$$p^{g_0}_{h_0}(n) \, 
= \, u(n)^{-1} \, (((t_1(n+g_0))^{10}(s_1(n+g_0))^{-9})^{10m_1(n)} \, ((t_2(n+g_0))^{10}(s_2(n+g_0))^{-9})^{-9m_2(n)})^{-9f_2(n+h_0)} \, =$$
$$= \, \hat u(n)^{-1} \, ((t_2(n+g_0))^{10}(s_2(n+g_0))^{-9})^{-9m_2(n)}$$ 
where the sequences, $\{t_j(n)\}$, $\{s_j(n)\}$, $j=1,2$, and $\{\hat u(n)\}$, are test sequences, and the
sequence $\{m_1(n),m_2(n)\}$ is a sequence of positive integers for which for every index $n$, 
$10m_1(n)-9m_2(n)=1$.

\medskip
At this stage we essentially repeat what we did so far with the   previous values of the pair $p,q$ (that were only in the set $E^1_1$ by
possibility (2)), for the new sequence of values of
this pair that are now known to be in the intersection $E^1_1 \cap E^2_1$.

By definition,
the group that is generated by the elements, 
$\{t_j,s_j\}_{j=1}^2$ and $\hat u,p,q$,  together with elements that are (formally) added for the powers:
$((t_j)^{10}(s_j)^{-9})^{10m_j}$, $j=1,2$,  
(each such formal power commutes with the corresponding element: 
$(t_j)^{10}(s_j)^{-9}$, $j=1,2$),  
is a completion, 
with a bottom level which is a free group of rank 5 (that is generated
by the elements,
$t_j,s_j$, $j=1,2$ and $\hat u$,  
and an upper level, in which  there are 2 abelian vertex groups of rank 2 that are
connected to a distinguished vertex.
We denote this
completion, $Comp_1$.

The subgroup  $Comp_1$ is a completion, and  the given sequence of values ($t_j(n+g_0),s_j(n+g_0)$, $j=1,2$, $ \hat u(n)$, $m_1(n),m_2(n)$)
is a test sequence of values of $Comp_1$, that projects to
values of the pair, $(p,q)$, that are in $NE(p,q)$, and  in $E^1_1(p,q) \cap E^2_1$. 
By the construction of envelopes of definable sets, as they
are presented in theorem 1.3 in  [Se6], with the intersection, $E^1_1 \cap E^2_1$, and the collection of all the test
sequences of $Comp_1$ that restrict to values of the pair $p,q$ that are in $E^1_1 \cap E^2_1$, 
it is possible to canonically associate finitely many envelopes,
that are closures of the completion, $Comp_1$. Since there exists a test sequence of $Comp_1$ that restricts to values of the pair
$p,q$ that are in the set $E^1_1 \cap E^2_1$, this collection of envelopes is not empty.

Hence, with  the  intersection $E^1_1 \cap E^2_1$, there are associated closures of
$Comp_1$, that we denote $Cl^1_1,\ldots,Cl^1_{d_1}$, that form the envelopes of the sets $E^1_1 \cap E^2_1$, with
respect to generic points (test sequences) of the completion, $Comp_1$. 
Furthermore, with each closure, $Cl^1_r$, there exists a finite (possibly empty)
closures of it, $Cl^1_{r,1},\ldots,Cl^1_{r,f_r}$. By the properties of envelopes (see theorem 1.3 in [Se6]),
the set of closures $Cl^1_1,\ldots,Cl^1_{d_1}$ and their associated closures,$Cl^1_{1,1},\ldots,Cl^1_{r,f_r}$, specifies precisely 
what test sequences of the completion, $Comp_1$, restrict to values of the variables $(p,q)$ that are in the set,
$E^1_1 \cap E^2_1$.

We now look at a subsequence of the original   test sequence, 
$t_j(n),s_j(n)$, $j=1,2$, and $\hat u(n)$,  
together with  a specific sequence of positive integers as values for the elements,
$\{m_1(n),m_2(n)\}$,
for which:  $10m_1-9m_2=1$. 

With the collection of all the test sequences of the completion, $Comp_1$, we have associated finitely
many envelope closures. 
With each of these closures 
there is an associated finite index subgroup of some fixed f.g.\ free abelian group, so we may look
at the intersection of all these finite index subgroups. Hence, with each closure there are
finitely many associated cosets of this finite index subgroup. By passing to a subsequence we may assume that
all the elements in the specific test sequence that we chose belong to a fixed coset of that finite index subgroup.

Now, for each index $n$ (of the specific test sequence that we chose), we look at the following sequence
of values for $1 \leq g,h \leq N+1$ (that is once again similar to the one that we used in proving lemma 5):
$$q^g_h(n) \, = \, ((t_1(n+g))^{10}(s_1(n+g))^{-9})^{10m_1(n+h)} \, \hat u(n)$$ 
$$p^g_h(n) \, = \, \hat u(n)^{-1} \, ((t_1(n+g))^{10}(s_1(n+g))^{-9})^{-9m_2(n+h)} .$$ 

By the way the values of the  sequence is defined,
and like in lemma 5, for every index $n$:
$(p^g_h(n),q^{g'}_{h'}(n)) \in NE(p,q)$, if and only if $g \neq g'$ or $g=g'$ and $h = h'$ 
($1 \leq g,g',h,h' \leq N+1$). Furthermore, 
since we have chosen the positive integers, $m_1(n),m_2(n)$, to be from a fixed 
class of the finite
index subgroup that is associated with the envelope closures of $E^1_2 \cap E^2_1$ with respect to the completion
$Comp_1$,
for every index $n$, and every $g \neq g'$, $1 \leq h,h',g,g' \leq N+1$, 
$(p^g_h(n),q^{g'}_{h'}(n)) \in E^1_1 \cap E^2_1$.

At this stage we apply the same argument that was used to prove lemma 5, 
to the 
equational set $E^1_1 \cap E^2_1$, and the co-equational set, $(E^1_2 \cap E^2_2)^c$, 
and obtain a contradiction. Hence, $NE(p,q)$ is not the union of two intersections between equational and co-equational sets,
and proposition 6 follows.

\line{\hss$\qed$}

Proposition 6 proves that the set $NE(p,q)$ is not the union of two intersections of equational and 
co-equational sets. Theorem 7 generalizes the proof to prove that $NE(p,q)$ is not a finite union of
such intersections.

\vglue 1pc
\proclaim{Theorem 7} The set $NE(p,q)$ is not  a finite union of intersections between an equational and 
a co-equational sets.
\endproclaim

\nfp The argument that we use is a straightforward generalization of the argument that was used to prove proposition 6.
Suppose that: $$NE(p,q) \, = \, (E^1_1 \, \cap \, (E^1_2)^c) \, \cup \, \ldots \, \cup \, 
(E^{\ell}_1 \, \cap \, (E^{\ell}_2)^c)$$ 
Let $M_{a_1,\ldots,a_b}(p,q) \, = \, E^{a_1}_1 \, \cap \ldots \cap \, E^{a_b}_1$, and 
    $T_{a_1,\ldots,a_b}(p,q) \, = \, E^{a_1}_2 \, \cap \ldots \cap \, E^{a_b}_2$,  
for $1 \leq b \leq \ell$,
and $1 \leq a_1 < \ldots < a_b \leq \ell$. Since the sets, $E^a_1$ and $E^a_2$, $a=1,\ldots,\ell$, are equational, so are the sets
$M_{a_1,\ldots,a_b}$  and
$T_{a_1,\ldots,a_b}$. 
Let $N$ be a bound on the lengths of 
strictly decreasing sequences: $\{\cap_{j=1}^f M_{a_1,\ldots,a_d}(p,q_j)\}_{f=1}^r$ and
 $\{\cap_{j=1}^f T_{a_1,\ldots,a_d}(p,q_j)\}_{f=1}^r$,
for all $1 \leq d \leq \ell$ 
and $1 \leq a_1 < \ldots < a_d \leq \ell$. 

We modify the argument that was used to prove proposition 6, and the definition of the elements $p^g_h$ and $q^g_h$ in the proof of that proposition.
We set: 
$$v^0_1=(t_1)^{10}(s_1)^{-9},\ldots,v^0_{2^{\ell}}=(t_{2^{\ell}})^{10}(s_{2^{\ell}})^{-9}$$
and then define the upper level elements iteratively:
$$v^1_1=(v^0_1)^{10m^1_1}(v^0_2)^{-9k^1_1},\ldots,
  v^1_{2^{\ell-1}}=(v^0_{2^{\ell}-1})^{10m^1_{2^{\ell-1}}}(v^0_{2^{\ell}})^{-9k^1_{2^{\ell-1}}}$$
$$v^r_1=(v^{r-1}_1)^{10m^r_1}(v^{r-1}_2)^{-9k^r_1},\ldots,
  v^r_{2^{\ell-r}}=(v^{r-1}_{2^{\ell-r+1}-1})^{10m^1_{2^{\ell-r}}}(v^{r-1}_{2^{\ell-r+1}})^{-9k^1_{2^{\ell-r}}}$$
for $1 \leq r \leq \ell$. Note that finally: 
$v^{\ell}_1=(v^{ell-1}_1)^{10m^{\ell}_1}(v^{\ell-1}_2)^{-9k^{\ell}_1}$.
We further set $q=(v^{\ell-1}_1)^{10m^{\ell}_1} \, u$ and 
$p  =  u^{-1} \, (v^{\ell-1}_2)^{-9k^{\ell}_1}$.  

We look at the collection of all the test sequences of values for the variables, 
$\{t_j,s_j\}_{j=1}^{2^{\ell}}$ and $u$.
Recall that elements in a test sequence need to satisfy the properties that are listed in definition 1.20 in [Se2],
they satisfy a very small cancellation assumption, and should be 
considered "generic" elements. For all the values in such test sequences, and arbitrary positive integer values
of the elements, $\{m^r_j\}$ and $k^r_j$, the corresponding values of the variables $(p,q)$ 
are in the set $NE(p,q)$.

By our assumptions, for each $(p_0,q_0) \in NE(p,q)$, there exists (at least one)  index $a$, $1 \leq a \leq \ell$,
for which $(p_0,q_0) \in \, E^a_1 \cap (E^a_2)^c$. Hence, in particular, $(p_0,q_0) \in E^a_1$. For each index
$a$, $1 \leq a \leq \ell$, we look at the collection of all the  test sequences of elements, 
$\{t_j,s_j\}_{j=1}^{2^{\ell}}$ and $u$,
and  sequences of positive integers,
$\{m^r_j\}$ and $\{k^r_j\}$,
for which each of the sequences, 
$\{m^r_j\}$ and $\{k^r_j\}$,
 converges to infinity,
and the corresponding sequence of values of the pair $(p,q)$ are all in the set $E^a_1$.

By the iterative definition of the elements $v^r_j$,
the group that is generated by the elements: $t_j,s_j$, $j=1,\ldots,2^{\ell}$, $u$, 
and $v^r_j$, $0 \leq r \leq \ell-1$, $1 \leq j \leq 2^{\ell-r}$, together with elements that represent
formal 
powers of the elements $v^r_j$ (each such formal power commutes with the corresponding element $v^r_j$),   
is a completion (i.e., an $\omega$-residually free tower - see section 6 in [Se1] and definition 1.12 in [Se2]).
This completion has a bottom level which is a free group of rank $2^{\ell+1}+1$ (that is generated
by $u$ and
$\{t_j,s_j\}_{j=1}^{2^{\ell}}$), and $\ell$ upper levels, where at level $r$, $1 \leq r \leq \ell$, there are $2^{\ell-r+1}$
abelian vertex groups, each of them of rank 2 and connected to the lower level with a cyclic edge group (that is generated by one of the
elements $v^r_j$). We denote this
completion, $Comp$.

The subgroup  $Comp$ is a completion, and the projection of generic points (test sequences) in its
associated variety project to values of the pair, $(p,q)$, that are in $NE(p,q)$, and hence are in $E^a_1(p,q)$  
for some index $a$, $1 \leq a \leq \ell$. Therefore, by the construction of envelopes of definable sets, as they
are presented in section 1 of [Se6] (see theorem 1.3 in [Se6]), with some of the groups, $E^a_1(p,q)$, 
there are associated envelopes,
that are closures of the completion, $Comp$.

Hence, with a (non-empty) subcollection of the  sets $E^1_1,\ldots,E^{\ell}_1$, there are associated closures of
$Comp$, that we denote $Cl_1,\ldots,Cl_d$ (which are the envelopes of the completion $Comp$ with respect to 
the collection of sets, $\{E^a_1\}$). Furthermore, with each closure, $Cl_g$, there exists a finite (possibly empty)
closures of it, $Cl_{g,1},\ldots,Cl_{g,f_g}$. By the properties of envelopes (see theorem 1.3 in [Se6]),
the set of closures $Cl_1,\ldots,Cl_d$ and their associated closures,$Cl_{1,1},\ldots,Cl_{g,f(g)}$, specifies precisely 
what test sequences of the completion, $Comp$, restrict to values of the variables $(p,q)$ that are in the sets,
$E^1_1,\ldots,E^{\ell}_1$.
% (and what test sequences restrict to values of $(p,q)$ that are
%in the complements of each of these sets). 

\medskip
We now look at a specific  test sequence of values for the variables, 
$\{t_j(n),s_j(n)\}_{j=1}^{2^{\ell}}$ and $u(n)$, 
together with  a specific sequences of positive integers,
$\{m^r_j(n)\}$ and $\{k^r_j(n)\}$, for which the integer values, $m^r_j(n)$ and $k^r_j(n)$, do not depend
on the indices $r$ or $j$,
$10m^r_j(n)-9k^r_j(n)=1$, and 
the sequences $m^r_j(n)$ and $k^r_j(n)$,
 converge to infinity. 

With the collection of all the test sequences of the completion, $Comp$, we have associated finitely
many envelope closures, $Cl_1,\ldots,Cl_d$, and with each closure $Cl_g$ we have further associated
finally many (possibly no) closures of it,  $Cl_{g,1},\ldots,Cl_{g,f(g)}$. With each of these closures 
there is an associated finite index subgroup of some fixed f.g.\ free abelian group, so we may look
at the intersection of all these finite index subgroups. Hence, with each closure there are
finitely many associated cosets of this finite index subgroup. By passing to a subsequence we may assume that
all the elements in the specific test sequence that we chose belong to a fixed coset of that finite index subgroup.

For each index $n$ (of the specific test sequence that we chose), we look at the following sequence
of values (that is similar to the one that we used in proving proposition 6 and lemma 5):
$$\{q^g_h(n)=(v^{\ell-1}_1(n+g))^{10m^{\ell}_1(n+h)} \, u(n)\}_{g,h=1}^{N+1} \ ; \
\{p^g_h(n)=u(n)^{-1} \, (v^{\ell-1}_1(n+g))^{-9k^{\ell}_1(n+h)}\}_{g,h=1}^{N+1}.$$  

By the way the elements $v^r_j$ are defined, and since we chose the elements 
$\{t_j(n),s_j(n)\}$, $m^r_j(n)$ and $k^r_j(n)$ ,from our fixed test sequence, like in proposition  6 and lemma 5:
$(p^g_h(n),q^{g'}_{h'}(n)) \in NE(p,q)$, if and only if $g \neq g'$ or $g=g'$ and $h = h'$. 

Since we have chosen the positive integers, $m^r_j(n)$ and $k^r_j(n)$, to be from a fixed class of the finite
index subgroup that is associated with the closures, $Cl_1,\ldots,cl_d$, and
$Cl_{1,1},\ldots,Cl_{d,f(d)}$, there exists a positive integer $b$, $1 \leq b \leq \ell$,
 and indices $1 \leq a_1 < \ldots < a_b \leq \ell$, so that for every index $n$, the pair $(q^g_h(n),p^{g'}_{h'}(n)) \in E^a_1$
for $g \neq g'$ and $1 \leq g,g',h,h' \leq N+1$, if and only if $a$ is one of the indices $a_1,\ldots,a_b$.

If $b=\ell$, then we apply the same argument that was used to prove lemma 5, 
to the 
equational set $E^1_1 \cap \ldots \cap E^{\ell}_1$, and the co-equational set, $(E^1_2 \cap \ldots \cap E^{\ell}_2)^c$, 
and obtain a contradiction. Hence, we may assume that $1 \leq b < \ell$. 

Let $M_b=
E^{a_1}_1 \cap \ldots \cap E^{a_b}_1$, and $T_b= 
E^{a_1}_2 \cap \ldots \cap E^{a_b}_2$.  
By the equationality of the set $T_b$, 
and the argument that was used to
prove lemma 5,  there  exists an index $g_0$, $1 \leq g_0 \leq N+1$, for which for every $h,h'$, $1 \leq h,h' \leq N+1$, 
$(p^{g_0}_h(n),q^{g_0}_{h'}(n)) \in T_b$. Since
$(p^{g_0}_h(n),q^{g_0}_{h'}(n)) \in NE(p,q)$ if and only if $h = h'$, the equationality of the set
$T_b$ and the argument that was used in proving lemma 5, imply that there exists an index $h_0$, 
$1 \leq h_0 \leq N+1$,  for which for every index $h'$, $1 \leq h' \leq N+1$,
$(p^{g_0}_{h'}(n),q^{g_0}_{h_0}(n)) \in T_b$. In particular,  
$(p^{g_0}_{h_0}(n),q^{g_0}_{h_0}(n)) \in T_b$. Since   
$(p^{g_0}_{h_0}(n),q^{g_0}_{h_0}(n)) \in NE(p,q) \cap M_b \cap T_b$  and:
$$NE(p,q) \, = \, (E^1_1 \, \cap \, (E^1_2)^c) \, \cup \ldots \cup \, 
(E^{\ell}_1 \, \cap \, (E^{\ell}_2)^c)$$ 
it finally follows that for each index $n$, there exists an index $a(n)$, $1 \leq a(n) \leq \ell$, and $a(n) \neq a_1,\ldots,a_b$ for which:
$(p^{g_0}_{h_0}(n),q^{g_0}_{h_0}(n)) \in M_b \cap E^{a(n)}_1$. By passing to a further subsequence of indices
(still denoted $n$),
we may clearly assume that $a(n)$ is independent of the index $n$, and that $a=a(n) \neq a_1,\ldots,a_b$.

Now, note that like what we did in proving proposition 6, the pairs,  
$(p^{g_0}_{h_0}(n),q^{g_0}_{h_0}(n))$,  can be written as: 
$$q^{g_0}_{h_0}(n)=(v^{\ell-1}_1(n+g_0))^{10m^{\ell}_1(n+h_0)} \, u(n) \, = \,
(v^{\ell-2}_1(n+g_0))^{10m^{\ell-1}_1(n+g_0)} \hat u(n) $$
$$p^{g_0}_{h_0}(n)=  u(n)^{-1} \, (v^{\ell-1}_1(n+g_0))^{-9k^{\ell}_1(n+h_0)} \, = \,  
\hat  u(n)^{-1} \, (v^{\ell-2}_2(n+g_0))^{-9k^{\ell-1}_1(n+g_0)}.$$

\medskip
Like the argument that was used in proving proposition 6, at this stage we essentially repeat what we did so far with the   previous values of the pair $p,q$ (that were only in the set $M_b$), 
for the new sequence of values of
this pair that are now known to be in the intersection $M_b \cap E^a_1$.

By the iterative definition of the elements $v^r_j$,
the group that is generated by the elements: $t_j,s_j$, $j=1,\ldots,2^{\ell-1}$, $u$, 
and $v^r_j$, $0 \leq r \leq \ell-2$, $1 \leq j \leq 2^{\ell-1-r}$, together with elements that represent
formal 
powers of the elements $v^r_j$ (each such formal power commutes with the corresponding element $v^r_j$),   
is a completion.
This completion has a bottom level which is a free group of rank $2^{\ell}+1$ (that is generated
by $u$ and
$\{t_j,s_j\}_{j=1}^{2^{\ell-1}}$), and $\ell-1$ upper levels, where at level $r$, $1 \leq r \leq \ell-1$, there are $2^{\ell-r}$
abelian vertex groups, each of them of rank 2 and connected to the lower level with a cyclic edge group (that is generated by one of the
elements $v^r_j$). We denote this
completion, $Comp_1$.

The subgroup  $Comp_1$ is a completion, and  the given sequence of values ($t_j(n+g_0),s_j(n+g_0)$, $j=1,\ldots,2^{\ell-1}$, $ \hat u(n)$, 
$v^r_j(n+g_0), m^r_j(n+g_0),k^r_j(n+g_0)$,  
$0 \leq r \leq \ell-2$, $1 \leq j \leq 2^{\ell-1-r}$) 
is a test sequence of values of $Comp_1$, that projects to
values of the pair, $(p,q)$, that are in $NE(p,q)$. Furthermore, by our conclusion, these  values of the pair $(p,q)$ 
are in the set $M_b \cap E^a_1$, where $1 \leq a \leq \ell$, and $a \neq a_1,\ldots,a_b$.

The values of the pair $(p,q)$ in our given sequence are contained in $M_b \cap E^a_1$. We pass to a further subsequence, so that
the values of $(p,q)$ are contained in a maximal intersection $M_b \cap E^{e_1}_1 \ldots \cap E^{e_s}_1$, where
$1 \leq e_1 <e_2 < \ldots <e_s \leq \ell$, $1 \leq s$, and each of elements $e_i \neq a_1,\ldots, a_b$. We set $b'=b+s$ and
$M_{b'}=M_b \cap E^{e_1}_1 \ldots \cap E^{e_s}_1$. Note that $1 \leq b < b' \leq \ell$.
Since we assume that the set, $e_1,\ldots,e_s$, is a maximal set, no proper subsequence of the given set of values of the
pair $(p,q)$, is contained in a set $M_{b'} \cap E^f_1$, where $1 \leq f \leq \ell$, and $f \neq a_1,\ldots,a_b,e_1,\ldots,e_s$.
 
By the construction of envelopes of definable sets, as they
are presented in theorem 1.3 in  [Se6], with the equational set $M_{b'}$,
and the collection of all the test
sequences of $Comp_1$ that restrict to values of the pair $p,q$ that are in $M_{b'}$, 
it is possible to canonically associate finitely many envelopes,
that are closures of the completion, $Comp_1$. Since there exists a test sequence of $Comp_1$ that restricts to values of the pair
$p,q$ that are in the set $M_{b'}$, this collection of envelopes is not empty.

Hence, with  the  equational set $M_{b'}$, there are associated closures of
$Comp_1$, that  form the envelopes of the set $M_{b'}$, with
respect to generic points (test sequences) of the completion, $Comp_1$. 
Furthermore, with each such envelope closure  there exists a finite (possibly empty)
collection of closures of it. By the properties of envelopes (see theorem 1.3 in [Se6]),
the set of envelope closures and their associated closures, specifies precisely 
what test sequences of the completion, $Comp_1$, restrict to values of the variables $(p,q)$ that are in the set,
$M_{b'}$.

We now look at  the original   test sequence of the completion $Comp_1$. 
With the collection of all the test sequences of the completion, $Comp_1$, we have associated finitely
many envelope closures. 
With each of these closures 
there is an associated finite index subgroup of some fixed f.g.\ free abelian group, so we may look
at the intersection of all these finite index subgroups. Hence, with each closure there are
finitely many associated cosets of this finite index subgroup. By passing to a subsequence of the original test
sequence of $Comp_1$, we may assume that
all the elements in the specific test sequence that we chose belong to a fixed coset of that finite index subgroup.

Now, for each index $n$ (of the specific test sequence that we chose), we look at the following sequence
of values for the pair $(p,q)$ and for indices, $1 \leq g,h \leq N+1$ 
(that is once again similar to the one that we used in proving lemma 5):
$$q^{g}_{h}(n) \, = \,
(v^{\ell-2}_1(n+g))^{10m^{\ell-1}_1(n+h)} \hat u(n) $$
$$p^{g}_{h}(n) \, =  \,   
\hat  u(n)^{-1} \, (v^{\ell-2}_1(n+g))^{-9k^{\ell-1}_1(n+h)}.$$

By the way the values of the  sequence is defined,
and like in proposition 6 and lemma 5, for every index $n$:
$(p^g_h(n),q^{g'}_{h'}(n)) \in NE(p,q)$, if and only if $g \neq g'$ or $g=g'$ and $h = h'$ 
($1 \leq g,g',h,h' \leq N+1$). Furthermore, 
since we have chosen the positive integers, $m^r_j(n),k^r_j(n)$, to be from a fixed 
class of the finite
index subgroup that is associated with the envelope closures of $M_{b'}$ with respect to the completion
$Comp_1$,
for every index $n$, and every $g \neq g'$, $1 \leq h,h',g,g' \leq N+1$, 
$(p^g_h(n),q^{g'}_{h'}(n)) \in M_{b'}$. 

$1 \leq b < b' \leq \ell$. If $b'=\ell$, we get a contradiction by the argument that was used to prove lemma 5. Hence, $b' < \ell$.
In that case we repeat the argument that we used for the first test sequence of values of $(p,q)$, and use the equationality of
$M_{b'}$ to prove the existence of indices, $g_0,h_0$ for every index $n$, so that  for every index $n$
there exists an index $a(n)$, $1 \leq a(n) \leq \ell$, and $a(n) \neq a_1,\ldots,a_b,e_1,\ldots,e_s$ for which:
$(p^{g_0}_{h_0}(n),q^{g_0}_{h_0}(n)) \in M_b \cap E^{a(n)}_1$. 

Therefore, after passing to a subsequence and repeating the same argument that we used for the first test sequence,
it is possible to enlarge the index $b'$ to $b''=b'+s'$, and obtain a completion $Comp_2$ with a test sequence that restricts
to
values of the pair $(p,q)$ that are in the set $M_{b''}$. Repeating this argument iteratively, after less than 
$\ell-1$ steps we can construct a sequence with similar properties, so that the values of the pair
$(p,q)$ is in the set $E^1_1 \cap \ldots \cap E^{\ell}_1$. Hence the same argument that was used to prove lemma 5, that
is applied to the equational set $E^1_1 \cap \ldots \cap E^{\ell}_1$ and the co-equational set 
$E^1_2 \cap \ldots \cap E^{\ell}_2$ gives a contradiction, and finally proves theorem 7.

\line{\hss$\qed$}

Since the union and the intersection of equational sets are equational, theorem 7 proves that $NE(p,q)$ is not in the
Boolean algebra of equational sets, hence that the theory of a free group is not equational. This concludes the proof of theorem
1.

\line{\hss$\qed$}

The results of [Se4] allow one to modify 
the argument that was used to prove theorem 1, and generalize theorem 1 to every non-elementary, torsion-free hyperbolic group.

\vglue 1pc
\proclaim{Theorem 8} The elementary theory of a non-elementary, torsion-free hyperbolic group is
not equational.
\endproclaim

\nfp To modify the argument that was used in the free group case, we only need to modify the proof of
theorem 7. By the results of [Se4], the same sieve procedure that was used for quantifier elimination
over a free group can be used to obtain quantifier elimination over a general torsion-free, non-elementary
hyperbolic group. By the results of [Se6], the sieve procedure and  the description of definable sets over 
a torsion-free hyperbolic group enable one to generalize the concept of an envelope of a definable set
from a free group to a general torsion-free, non-elementary hyperbolic group.

As there exist quasi-isometric embeddings of a non-abelian free group into any given non-elementary torsion-free 
hyperbolic group, test sequences  generalize to torsion free hyperbolic groups (and indeed, test sequences in
hyperbolic groups are used extensively in [Se4]).
Therefore, to generalize the argument that was used to prove theorem 7 to non-elementary, torsion-free
hyperbolic groups, the only result that still needs to be generalized to hyperbolic groups is
the Lyndon-Schutzenberger theorem [Ly-Sch], that shows that for $m,n,p \geq 2$, the solutions
of the equation $x^my^nz^p=1$ in a free group, generate a cyclic subgroup. For the purposes of generalizing 
the proof of theorem 7 we need a slightly weaker result that can be proved easily using Gromov-Hausdorff
convergence.

\vglue 1pc
\proclaim{Lemma 9} Let $\Gamma$ be a non-elementary, torsion-free hyperbolic group. There exist 
some positive integers $\ell,b>0$, so that for every $m>b$, the solutions of the system $x^{\ell+1}y^{-\ell}z^m=1$ in $\Gamma$,
generate a cyclic subgroup in $\Gamma$.
\endproclaim

\nfp Let $\Gamma$ be a non-elementary, torsion-free hyperbolic group for which such positive integers $\ell$ and $b$ do not
exist. In that case there exist  triples $x_n,y_n,z_n \in \Gamma$, for an increasing sequence of integers $n$, and for which: 
\roster
\item"{(1)}" $x_n^{n+1}y_n^{-n}z_n^m=1$ in $\Gamma$.

\item"{(2)}" $m>n$, and $x_n,y_n,z_n$ do not belong to a cyclic subgroup in $\Gamma$, hence, they do not commute.
\endroster
The sequence of values $x_n,y_n,z_n$, define a sequence of homomorphisms from a free group $F$ generated by  $<x,y,z,u,v,w>$, to $\Gamma$,
$h_n:F \to \Gamma$, $h_n(x)=x_n$, $h_n(y)=y_n$, $h_n(z)=z_n$, $h_n(u)=x_n^{n+1}$, $h_n(v)=y^{-n}$, $h_n(w)=z_n^m$. This
sequence of actions subconverges into a non-trivial  action (not necessarily faithful) of $F$ on some real tree $Y$. 
 $x_n^{n+1}y_n^{-n}z_n^m=1$ in $\Gamma$, and this equation corresponds to a (thin) triangle that is associated with the equation:
$h_n(u)h_n(v)h_n(w)=1$ in $\Gamma$. In the limit tree $Y$, this sequence of thin triangles converges in the Gromov-Hausdorff topology
to a (possibly degenerated) tripod, and this tripod contains a (non-degenerate) segment that is stabilized by at least two elements
from the triple $x_n,y_n,z_n$. By [Pa], this implies that for large enough index $n$ (from the convergent subsequence), two of the elements
$x_n,y_n,z_n$ commute, hence, these two elements are contained in the same cyclic subgroup of $\Gamma$. In a torsion-free hyperbolic group, 
this implies that for large enough $n$, all the 3 elements, $x_n,y_n,z_n$, belong to the same cyclic subgroup of $\Gamma$, a contradiction to
(2). 

\line{\hss$\qed$}

Lemma 9 and the results of [Se4] and [Se6] allow one to generalize the proof of theorem 7 to non-elementary,
torsion-free hyperbolic groups, hence, to prove that such groups are not equational.

\line{\hss$\qed$}

%\newpage

\vskip 2em

\end